\documentclass[10pt]{article}
\usepackage{amsmath, amssymb}
\usepackage[usenames, dvipsnames]{color}
\begin{document}

\title{Correcting the proof of Theorem 3.2 
in {\it Almost sure rates of mixing for i.i.d. unimodal maps} 
(V. Baladi, M. Benedicks, V. Maume-Deschamps)  
Ann. E.N.S. (2002)}
\date{August 27, 2010}
\maketitle

 We thank Weixiao Shen, who pointed out to us  that
 the proof of Theorem ~3.2 of 
 \cite{BBM}
  was flawed: Since
 $(F^R_\omega)_* (m_0| \Delta_{\sigma^{-n}\omega,0}, n \le 1)$ 
 is not a probability measure on
 $\Delta_{\omega, 0}$
 in general, the claim three lines
 below (3.9) that
 the sequence $\phi_{n, \omega}$ 
 is \footnote{Note also the -- unimportant ---
 typo in the def. of  $\phi_{n, \omega}$ where
 $\frac{1}{n}$ should read $\frac {d}{dm_0}\frac{1}{n}$.}
  bounded, uniformly in $n$,
 is unfounded. (Indeed, even if  $F^R_\omega$ is supposed
 piecewise affine, counter-examples may be constructed.)
 
 Weixiao Shen kindly  provided the following argument  below to fix this proof.
 
 {\it  \color{RubineRed} We do not claim anymore  that $\{h_\omega^{-1}\} \in \cal F^{\cal K_\omega}_\beta$ if (3.2)
 holds.}
 See below how to deduce Lemma~5.1 from mixing of $F$ and
 upper bounds from $h_\omega$, without 
 using lower bounds for $h_\omega$,
 on which the proof of (5.1) depended. The lower bound for $h_\omega$
 is not used elsewhere.

\smallskip

{\bf Corrected proof of Thm 3.2:}

Instead of working
with $F^R_\omega$, we directly work with $F_\omega$
(just like in the proof of Sublemma 5.5 (1)). 
More precisely, for each $\omega$ and $n\ge 0$, let $\mu_n^\omega$ 
be the push-forward of
$m_0|\Delta_{\sigma^{-n}\omega,0}$ by 
$F^n_{\sigma^{-n}\omega}$.
This push-forward is a probability measure on the tower 
$\Delta_\omega$, absolutely continuous with
respect to $m$. Estimate (3.9) implies that the densities $\varphi_{n}^{\omega}$
of the
$\mu_n^\omega$ belong to $\cal F^+_\beta$, with 
constants $\sup_n C_{\varphi_{n,\omega}}<\infty$. 

Recall that (3.3) says that
$m(\Delta_\omega)<\infty$ for almost every
$\omega$.

In the application to unimodal maps, (3.3)
means that we may view almost every  $\Delta_\omega$
as a compact interval. Thus, Arzela-Ascoli
gives for almost every $\omega$
a subsequence $n_\ell \to \infty$
so that $\frac{1}{n_\ell}\sum_{k=0}^{n_\ell-1}\varphi_{k}^{\omega}$
converges to the density of
a probability measure on $\Delta_\omega$. 

In the general case,  (3.3) implies that for almost all
 $\omega$, there is
a subsequence $n_\ell \to \infty$
so that $\frac{1}{n_\ell}\sum_{k=0}^{n_\ell-1}\mu_{k}^{\omega}$
converges in the weak-(*) topology to
a probability measure on $\Delta_\omega$, absolutely continuous with respect
to $m$.

In both cases,
by the diagonal principle, for almost every $\omega$, we can find a sequence $n_m$ such that for each integer $N$
$$\frac{1}{n_m}\sum_{k=N}^{n_m-1} \mu_{k-N}^{\sigma^{-N}\omega}$$
converges to 
a probability measure $\mu_{\sigma^{-N}\omega}$
on the tower $\Delta_{\sigma^{-N}\omega}$, absolutely
continuous with respect to $m$. By construction, 
$(F_\omega)_*(\mu_\omega)=\mu_{\sigma \omega}$ for almost every
$\omega$. 
This 
gives the claimed 
absolutely continuous
sample sationary probability measure.

Next, the construction implies that the density $h_\omega$
of $\mu_\omega$ is bounded uniformly for almost all $\omega$.
In particular,
$\{h_\omega, \, \mbox{a.a.} \, \, \omega \}$  belongs to $\cal F^+_\beta
\cap \cal F^{\cal K_\omega}_\beta$, where $\cal K_\omega$ is
bounded uniformly over almost all $\omega$.

This ends the proof of Theorem 3.2, with the statement on $h_\omega^{-1}$
removed.
(The proof  does not require (3.2).)

 \medskip
For the sake of comparison with the original proof of Theorem 3.2, we note that
the restrictions $\hat \nu_\omega$
of $\mu_\omega$ to $\Delta_{\omega,0}$
give a family of finite measures, which is
invariant under $F^R_\omega$, in the sense that for almost all $\omega$
and each $E \subset \Delta_{\omega, 0}$
$$
\hat \nu_\omega(E)=
\sum_{\ell=1}^{\infty}  \hat \nu_{\sigma^{-\ell}\omega}
((F^R)^{-1}(E)\cap \Delta_{\sigma^{-\ell}\omega, 0})\, .
$$ 
(A priori  $\hat \nu_{\sigma^{-n}\omega}(\Delta_{\sigma^{-n}(\omega),0})=
\mu_{\sigma^{-n}(\omega)}(\Delta_{\sigma^{-n}(\omega),0})$ may depend on 
$\omega$
and $n$, i.e., our assumptions do not guarantee a common normalisation
factor.)

Note also that $h_\omega$ may vanish
at $(x,\ell)$ for $\ell \ge 1$, but then it vanishes identically
on the element of $\cal Z_\omega$ containing $(x,\ell)$.

\bigskip

{\bf Corrected proof of Lemma 5.1:}

First note that, instead of (5.1), we may use the weaker claim 
$\int_\Omega V_\omega^\ell dP>0$ to show that
$\int \exp [-\upsilon V^{\tau_i-\tau_{i-1}}_{\sigma^{\tau_{i-1}}}]
\, dP(\omega) <1$ if $\upsilon >0$
is  small enough. We explain how to deduce 
$\int_\Omega V_\omega^\ell dP>0$
from mixing of $F$:

Let \footnote{ $F^{-\ell}_\omega(\Delta_{\omega,0})$
should be replaced by  $F^{-\ell}_\omega(\Delta_{\sigma^{\ell}\omega,0})$
 in the
definition of $V_\omega^\ell$ on page 92 of the paper.} 
$\hat{V}_\omega^\ell =\mu_\omega(\Delta_{\omega,0}\cap
F^{-n}_\omega (\Delta_{\sigma^n\omega,0}))
$.
Then a simple application of
the mixing property of $F$ implies 
$$
\lim_{\ell \to \infty}
\int_\Omega \hat{V}_\omega^\ell dP(\omega)= \mu_{\epsilon}(\Lambda) ^2
\, .
$$
It is easy to prove that
$\mu_\epsilon(\Lambda)>0$.  
So we obtain a positive lower bound for 
$\int_\Omega \hat{V}_\omega^\ell dP(\omega)$,
for $\ell$ large. Since the density $h_\omega$
of $\mu_\omega$ is bounded from above,
$V_\omega^\ell > C \hat{V}_\omega^\ell$, so $\int V_\omega^\ell dP$ 
is bounded away
from zero.  This ends the proof.

\end{document}


\title{Correcting the proof of Corollary 5.2
in {\it Almost sure rates of mixing for i.i.d. unimodal maps} 
(V. Baladi, M. Benedicks, V. Maume-Deschamps)  
Ann. E.N.S. (2002)}
\date{July 1, 2017}
\maketitle

 We thank Wael Bahsoun, Christopher Bose, and Marks Ruziboev, who pointed out to us  that
 the proof of Corollary~5.2 in
 \cite[p. 94]{BBM}
  was flawed: In the inequality in the 9th line of this proof, 
  the supremum over $\omega$ of the sum must be replaced by the sum of the suprema
  over $\omega$. The next inequality then does not hold. We explain below how to
 recover Corollary 5.2
    {\it  \color{RubineRed} up to replacing the set $M'_q$ on p. 93  by
  the smaller set
  \begin{align*}
  M'_q=\bigl \{
  (\omega, x, x')\in \bigcup_{\omega \in \Omega} (\{\omega\} \times \Delta_\omega
  \times \Delta_\omega)
  \mid&
  \omega \in M_q^{\{ \tau_i^\omega(x,x'), 0\le i \le q\}} \\
 &\quad \mbox{ and } \tau_{q}^\omega(x,x') > [q^{1/v}]
  \bigr \}\, ,
  \end{align*}
  where the constant $v \in (0,1/4)$ is determined by 
   the estimate
  (I) in the proof of Proposition~5.6.}
  Since Corollary 5.2 is only used to obtain the estimate
  (II) in the proof of Proposition 5.6, and since $q\le n^v<q+1$ there (for the
  same $v$), the rest of the paper
  is not affected.

\smallskip

{\bf Corrected proof of Corollary 5.2:}
First note that Corollary 5.2 is not used to prove Lemma 5.3,
 Lemma 5.4, Sublemma 5.5, or the bound for (I)
in Proposition 5.6. We can thus postpone the statement
and proof of Corollary~5.2 to just  before the proof of Proposition 5.6.

Next, for  fixed $\omega$, $q$, and $\{\tau_i, i=1, \ldots , q\}$ with $\tau_{j}-\tau_{j-1}\ge \ell_0$,
set $$(\Delta_\omega^2)^{\{\tau_i\}}=
(\Delta_\omega \times \Delta_\omega) \cap \{\tau_i^\omega(x,x')=\tau_i, i=1, ... q\}\, .
$$

By adapting the argument used to bound (I) in the proof of Proposition 5.6
on pp. 97-98, using $n_5$ defined there,  we can prove that
 for any $\omega$, any $n>n_5(\omega)$, and any fixed $\tau_0, \tau_1,\ldots ,\tau_q$
 so that  $\tau_{q-1}\le n <\tau_q$, we have
(set $\tau_i-\tau_{i-1}=k_i$)
$$m\times m ((\Delta_\omega^2)^{\{\tau_i\}}) < C(\epsilon) e^{-n^1/4}\, .$$
Summing over the $k_i$ so that $\sum_{i=1}^q k_i \le n$ like on p. 98, we still find
$C(\epsilon) e^{-n^{1/4}}$, up to changing the constant.
(This step is the key estimate.)

Since we can assume that ${\cal K}_\omega^2 < \kappa^{-q/8}$ (because
$P( {\cal K}_\omega^2 > \kappa^{-q/8})$ is small by (3.8)),
it suffices to  bound the sum of
$$ \kappa^q \kappa^{-q/8}
     \sum_{n> \tilde v^8 |\log \kappa|^8} \, \, \sum_{\tau_i: \tau_{q-1}\le n <\tau_q}\, \, 
      \sup_{\omega \in \Omega^{\{\tau_i\}}, n_5(\omega)<\tilde v^8 |\log \kappa|^8}
                            m\times m ((\Delta_\omega^2)^{\{\tau_i\}})
$$
(for which we can use the key estimate)
 with
$$\kappa^q \kappa^{-q/8}
  \sum_{n< \tilde v^8 |\log \kappa|^8} \, \, \sum_{\tau_i: \tau_{q-1}\le n <\tau_q}\, \, 
     \sup_{\omega \in \Omega^{\{\tau_i\}}} m\times m ((\Delta_\omega^2)^{\{\tau_i\}})
     $$
 (which is tractable by hand)
 and
   $P(n_5(\omega)> \tilde v^8 |\log \kappa|^8$) (which is $\le \kappa^{\tilde v}$ by p.~97).
We conclude by choosing $\tilde v >0$  suitably.